\documentclass[12pt, reqno]{amsart}

\usepackage[top=1in, bottom=1in, left=1.2in, right=1.2in]{geometry}

\usepackage{amsmath, amsthm, amssymb}

\usepackage{enumerate}
\usepackage{url}

\usepackage{fancyvrb}

\usepackage{booktabs}
\usepackage[table]{xcolor}
\usepackage{caption}
\usepackage{longtable}
\usepackage{hyperref}

\newtheorem{thm}{Theorem}[section]
\newtheorem{lem}[thm]{Lemma}
\newtheorem{prop}[thm]{Proposition}
\newtheorem{cor}[thm]{Corollary}

\theoremstyle{definition}

\theoremstyle{remark}

\numberwithin{equation}{section}

\newcommand{\mmod}[1]{\,\,\textnormal{mod}\,\,#1}

\def\alp{{\alpha}} 
\def\bet{{\beta}} 
\def\gam{{\gamma}}

\def\kap{{\kappa}}
\def\lam{{\lambda}}

\def\ome{{\omega}}  
\def\eps{\varepsilon}
\def\epsilon{\eps}
\def\implies{\Rightarrow}

\def \le {\leqslant} \def\ge{\geqslant}

\def\d{{\,{\rm d}}}

\def \bN {\mathbb N}

\def \bR {\mathbb R}
\def \bZ {\mathbb Z}
\def \bT {\mathbb T}

\def \bc {\mathbf c}

\def \bv {\mathbf v}

\def \fm {\mathfrak m}

\def \fJ {\mathfrak J}
\def \fK {\mathfrak K}

\def \fM {\mathfrak M}

\def \fS {\mathfrak S}

\def \cA {\mathcal A}

\def \cD {\mathcal D}

\def \cJ {\mathcal J}

\def \cP {\mathcal P}

\def \meas {\mathrm{meas}}

\renewcommand{\geq}{\ge}
\renewcommand{\leq}{\le}

\begin{document}
\title{A linear bound in Waring's problem}



\author{Sam Chow \and Rainer Dietmann}

\address{Mathematics Institute, Zeeman Building, University of Warwick, Coventry CV4 7AL, United Kingdom}
\email{Sam.Chow@warwick.ac.uk}

\address{Department of Computing, Security and Mathematics, Royal Holloway, University of London\\
Egham TW20 0EX, United Kingdom}
\email{Rainer.Dietmann@rhul.ac.uk}

\begin{abstract} 
We show that $G(k) \asymp k$ in Waring's problem.
\end{abstract}

\maketitle

\section{Introduction}

Waring's problem \cite{Waring1770}, \cite[\S1.1]{Vaughan97} dates back to 1770. In its modern formulation, it concerns $G(k)$, which is the least $s$ such that any sufficiently large positive integer is a sum of at most $s$ many perfect $k^{\mathrm{th}}$ powers.  A simple counting argument implies that 
\[
G(k) \ge k+1
\]
for all $k \ge 2$, see for example \cite[\S 21.6]{HW2008}. 
Pioneering work of Hardy and Littlewood \cite{HardyLittlewood25}
delivered a bound of the form $G(k) \ll k 2^k$, and Vinogradov \cite{Vinogradov35} later lowered the order of magnitude to $G(k) \ll k \log k$.
For large $k$, the best upper bound hitherto obtained is due to Br\"udern and Wooley \cite[Theorem 1.1]{BW2023}, who showed that
\[
G(k) \le \lceil k(\log k + 4.20032) \rceil
\]
for all $k \in \bN$.
This used the circle method of Hardy and Littlewood \cite{Vaughan97}, building on work of Davenport \cite{Davenport39, Davenport42}, \cite[\S 11]{VW2002}, and subsequent work of Vaughan \cite{Vau1989} and Wooley~\cite{Woo1992} involving smooth numbers.

We determine that $G(k)$ has order of magnitude $k$.

\begin{thm} 
\label{MainThm}
We have 
\[
G(k) \ll k
\]
for all $k \in \bN$.
\end{thm}

More sophisticated reasoning reveals that $G(k) \le 8.4 k$ for all $k \in \bN$.
To facilitate readability, we do not demonstrate this at present.

\subsection*{Methods}

We use the circle method, noting by orthogonality that the number of representations of $N$ as a sum of $s$ many positive $k^{\mathrm{th}}$ powers is
\[
R_{s,k}(N) = \int_\bT
\left( \sum_{n \le P} e(\alp n^k) \right)^s e(-N \alp) \d \alp,
\]
where $P = N^{1/k}$ and $\bT$ is the torus $\bR/\bZ$.  We will establish the stronger result that
\[
R_{s,k}(N) \gg_{s,k} P^{s-k}.
\]
The following mean value estimate is key.

\begin{thm} 
\label{MeanValue}
There exist an integer $r \ll k$, a constant $d_k > 0$, and $\cA_k(P) \subseteq [1,P] \cap \bZ$ such that
\[
\# \cA_k(P) \sim d_k P \quad (P \to \infty),
\qquad
\int_\bT \left|
\sum_{n \in \cA_k(P)} e(\alp n^k) \right|^r \d \alp \ll_k P^{r - k + 4/k^2}.
\]
\end{thm}

Additional smooth variables then deliver a further saving of more than $P^{4/k^2}$ on the minor arcs. The major arcs then provide the main contribution and complete the proof.

To establish Theorem \ref{MeanValue}, we use `factors and fibres'. The former is a $p$-adic analogue of the diminishing ranges argument, as espoused by Davenport \cite[\S 11]{VW2002}. Factoring out some $m^k$ from the phase $\alp n^k$ of a Weyl sum of length $X$, the length shortens from $X$ to $X/m$. Here $m$ is a product of primes from disjoint intervals, and $n/m$ is a smooth cofactor. This is a deliberate restriction on the variables, and enables smooth Weyl sums to enter the fray.

The latter refers to fibres of the multiplication by $m^k$ map on $\bT$. Assume that $m$ and $W$ are not too large. The main proposition then asserts that, for any $\bet \in \bR$, there are at most $O_k(W^2)$ many $\alp \in \bT$
such that 
\[
m^k \alp = \bet
\]
and the Weyl sum has size at least $X/W$. Its proof uses the Gram matrix of the phase vectors $e(n^k \alp)$ on the fibres, and we use congruences to control the entries.

The main proposition enables an iterative process on the large spectrum. It parallels Vaughan's \cite{Vau1989}, but is measure theoretic. Another difference is that our process is adaptive. It is important to carefully choose how many primes to remove at each step. The trade-off is between the effective shortening of the Weyl sums and the bound on the number of bad fibres.

\subsection*{Notation}
The letter $p$ is reserved for primes. For $x \in \bR$, we write $e(x) = e^{2 \pi i x}$.
We adopt the Bachmann--Landau and Vinogradov asymptotic notations.

\subsection*{AI statement}
The work was carried out by GPT-6.2 Astra,  rewritten by the authors, and refereed by GPT-6.2 Astra.

\subsection*{Acknowledgements}
We thank J\"org Br\"udern and Trevor Wooley for teaching us about Waring's problem and the circle method.

\subsection*{Funding}
SC is supported by a grant of use from OpenAI as part of the `ChatGPT for Academic Researchers' programme.

\subsection*{Rights}
For the purpose of open access, the authors have applied a Creative Commons Attribution (CC-BY) licence to any Author Accepted Manuscript version arising from this submission.

\section{Setup}
\label{setup}

We assume as we may that $k$ is arbitrarily large, and that $P$ is arbitrary large in terms of $k$. Put
\[
\gam = k^{-3}, \quad
L = \log(1/\gam) = 3 \log k, \quad
R = k^{10}, \quad
I = [1,R] \cap \bZ.
\]
For $i \in I$, let
\[
a_i = \frac{1 - 2 \gam}{R}
\left( 1 + \frac{i-(R+1)/2}{2R}
\right),
\]
so that
\[
a_i \asymp 1/R, \qquad
\sum_i a_i = 1 - 2 \gam.
\]
Note that, for some $\ome = \ome(k) \in (0,\gam/10)$,
\begin{equation}
\label{omega}
(1+\ome)a_i < a_{i+1} \qquad (1 \le i < R).
\end{equation}
For each $i \in I$, denote 
\[
\cP_i = \{ p: P^{a_i} < p \le P^{(1+\ome)a_i} \}.
\]
For $T \subseteq I$, write 
\[
\cD_T = \left\{ \prod_{i \in T} p_i : p_i \in \cP_i \right\}.
\]
Fix a sufficiently small positive constant $\eta < \min \{ a_i: i \in I \}$, and define
\[
f_Y(\alp) =
\sum_{\substack{
u \le Y \\
p \mid u \implies p \le Y^\eta}} e(\alp u^k).
\]

Let $S \subseteq I$,
and suppose $p_i \in \cP_i$ is already chosen for each $i \notin S$. Then define
\begin{equation}
\label{state}
X = \frac{P}{\prod_{i \notin S} p_i} = P^x,
\qquad
F_S(\alp; X) = \sum_{m \in \cD_S} f_{X/m} (m^k \alp), \qquad
F(\alp) = F_I(\alp; P).
\end{equation}
Note that $F$ is a trigonometric polynomial supported on the $k^{\mathrm{th}}$ powers of elements of
\[
\cA_k(P) = \{ mu \le P:
m \in \cD_I, \ u \text{ is } (P/m)^\eta\text{-smooth} \}.
\]

 \begin{lem} 
\label{combinatorics}
We have
\begin{equation}
\label{xineqs}
\gam \le x \le 1, \qquad
0 \le x - \sum_{i \in S} a_i \le 2 \gam.
\end{equation}
If $\emptyset \ne T \subseteq S$ then
\begin{equation}
\label{subset}
F_S(\alp; X) = \sum_{m \in \cD_T} F_{S \setminus T} (m^k \alp; X/m), \qquad
\sum_{m \in \cD_T} \frac1m \le 1.
\end{equation}
For $T_1 \subseteq T$ and $m_1 \in \cD_{T_1}$,
\[
\sum_{\substack{ m \in \cD_T \\
m_1 \mid m }} F_{S \setminus T} (m^k \alp; X/m) = F_{S \setminus T_1} (m_1^k \alp; X/m_1).
\]
\end{lem}

\begin{proof} The inequalities \eqref{xineqs} follow from the observation that
\[
x - \sum_{i \in S} a_i = 2 \gam - \sum_{i \notin S} (\log_P p_i - a_i) \in [2 \gam - \ome, 2 \gam].
\]

The smooth cofactor has no divisors in any $\cP_i$, so any integer $n = mu$ in the sum $F_S(\alp; X)$ determines its factorisation. The identity in \eqref{subset} follows by grouping according to the primes indexed by $T$, and the final identity follows similarly.

By Mertens's theorem,
\[
\sum_{p \in \cP_i} p^{-1} \sim \log (1 + \ome) \qquad (P \to \infty),
\]
Taking the product over $T$, each factor is $1$-bounded, so the product is too. This confirms the inequality in \eqref{subset}. 
\end{proof}

\section{Large values on a fibre}

In this section, we will establish our main proposition. For the whole section, let us assume
that
\[
P^{k^{-3}} = P^\gamma \le X \le P, \qquad 1 \le M \le L=3 \log k,
\]
and that
\[
  m \le X^{M/k}
\]
is a positive, square-free integer that is a product of at most $R=k^{10}$ many primes exceeding $k$.
Note that for fixed degree $k$, for sufficiently large $P$ all primes in the sets $\cP_i \; (i \in I)$ satisfy this property of exceeding $k$.

Our next lemma is a consequence of a Weyl-type inequality.

\begin{lem}
\label{RationalFrequency}
Put $c = 1 / 1000$. Suppose $q \mid m^k$ and $q>X$. Further, let $a$ be an integer coprime to $q$. Then
\begin{equation}
\label{weyl}
\sum_{\substack{n \le X \\ (n,m) = 1}} e(an^k/q) \ll_k X^{1-c/M^2}.
\end{equation}
\end{lem}

\begin{proof}
Define
\[
  U=\log_X q, \qquad
  D = \lceil 4 U \rceil,
\]
and note that $1<U \le M.$ 
For a prime $p \mid m$, let
\[
  e_p=v_p(q)
\]
be the highest power of $p$ dividing $q$, and write
\[
  s_p=\max\left\{1,\left\lceil\frac{e_p}{D+1/2}\right\rceil\right\},
\qquad
  A = \prod_{p \mid m} p^{s_p}.
\]
We split the sum over $n$ into residue classes modulo $A$. As $m$ is
square-free and $s_p \ge 1$ for all $p \mid m$, we have $m \mid A$. Hence by the
summation condition $(n,m)=1$, for a residue classes of the shape
$n=n_0+Ay$, the condition $(n,m)=1$ is satisfied if and only of $(n_0,m)=1$. Therefore
\begin{equation}
\label{classes}
\sum_{\substack{n \le X \\ (n,m) = 1}} e(an^k/q) =
\sum_{\substack{1\leq n_0\leq A\\(n_0,m)=1}}
\ \sum_{\substack{y\geq 0\\n_0+Ay\leq X}}
e\left(\frac{a(n_0+Ay)^k}{q}\right).
\end{equation}

Now fix $n_0$ with $(n_0,A)=1$. The binomial theorem yields
\[
\frac{a(n_0+Ay)^k}{q}
=
\sum_{j=0}^k
\frac{a}{q}\binom{k}{j}
n_0^{k-j}A^j y^j.
\]
Observe that for $j>D$ the coefficient
\[
  \frac{a}{q}\binom{k}{j}
n_0^{k-j}A^j 
\]
is an integer. Indeed, suppose $p \mid q$. Then $p$ also divides $m$, and so
\[
  \nu_p(A^j)=j s_p \ge (D+1) s_p \ge (D+1) \frac{e_p}{D+1/2} \ge e_p.
\]

In \eqref{classes}, we may therefore replace $a(N_0+Ay)^k$ by its binomial
expansion stopped at power $y^D$. The coefficient of $y^D$ is
\[
\lambda_D
:=
\frac{a}{q}\binom{k}{D}
n_0^{k-D}A^D.
\]
We now determine the denominator of this
coefficient in lowest terms. Using 
\[
q \le m^k \le X^M \le X^L,
\]
we obtain
$
  D \le 4U + 1 \le 12 \log k+1 \le k,
$
as $k$ is large. Every prime divisor of $q$ divides $m$ and is
greater than $k$, and therefore does not divide
$\binom{k}{D}$.
Moreover
\[
(a,q)=1,\qquad (n_0,q)=1,
\]
the second assertion following from $(n_0,m)=1$
and $q$ dividing $m^k$. Consequently, cancellation in $\lambda_D$ can
arise only from the factor $A^D$.
Since
$
v_p(A^D)=Ds_p,
$
the denominator of $\lambda_D$, in lowest terms,
is
\[
q_D := \prod_{p\mid m}
p^{\max\{e_p-Ds_p,0\}}.
\]
This is the denominator to which we shall apply
a degree-$D$ Weyl-type estimate.

The definition of $s_p$ ensures that
\[
\frac{e_p}{D+1/2}
\leq s_p
\leq \frac{e_p}{D+1/2}+1.
\]
%The upper bound remains valid when $e_p=0$.
Taking logarithms in the definition of $A$, and using that $m$ is square-free, we
obtain
\[
\begin{aligned}
\log A
&=\sum_{p\mid m}s_p\log p \leq
\frac{1}{D+1/2}\sum_{p\mid m}e_p\log p
+\sum_{p\mid m}\log p=\frac{\log q}{D+1/2}+\log m.
\end{aligned}
\]
Dividing by $\log X$ yields
\begin{equation}
\label{logXA}
\log_X A
\leq\frac{U}{D+1/2}+\log_X m.
\end{equation}

Next, since
$
s_p\geq \frac{e_p}{D+1/2},
$
we have
\[
e_p-Ds_p
\leq e_p\left(1-\frac{D}{D+1/2}\right)
=\frac{e_p}{2D+1}
\]
and hence $
\max\{e_p-Ds_p,0\}
\leq\frac{e_p}{2D+1}.
$
Summing over the prime divisors of $m$ gives
\[
\log q_D\leq\frac{\log q}{2D+1}.
\]
For a lower on $\log q_D$, we use
$
s_p\leq\frac{e_p}{D+1/2}+1
$
to find
\[
\max\{e_p-Ds_p,0\} \ge e_p-Ds_p
\geq \frac{e_p}{2D+1}-D.
\]
Consequently,
\[
\log q_D
\geq
\sum_{p\mid m}
\left(\frac{e_p}{2D+1}-D\right)\log p=\frac{\log q}{2D+1}-D\log m.
\]
The upshot is that
\begin{equation}
\frac{U}{2D+1}-D\log_X m
\leq\log_X q_D
\leq\frac{U}{2D+1}.
\label{logXqD}
\end{equation}
We now make these estimates explicit.
As $U > 1$ and
$4U\leq D<4U+1,
$
\[
\frac1{11} < \frac{U}{8U+3}
<
\frac{U}{2D+1}
\leq\frac{U}{8U+1} <\frac18.
\]
On the other hand,
\[
D\log_X m
\leq (4L+1)\frac{M}{k}
\leq\frac{(4L+1)L}{k}
=o(1).
\]
As $k$ is large, it now follows from \eqref{logXqD} that
\begin{equation}
X^{1/12}\leq q_D\leq X^{1/7}.
\label{qDbound}
\end{equation}

Similarly, using $D\geq4U$ in \eqref{logXA}, we obtain
\[
\log_X A
\leq\frac{U}{D+1/2}+\frac{M}{k}
\leq\frac14+\frac{L}{k} \le \frac13
\]
and so $A \le X^{1/3}$. Now fix $n_0 \in [1,A]$ coprime to $m$ and let
\begin{equation}
Y=Y_{n_0}
=\#\{y\geq0:n_0+Ay\leq X\}
\asymp\frac{X}{A}
\gg X^{2/3}.
\label{Ylower}
\end{equation}

By the degree-$D$ Weyl-type
estimate of Bourgain~\cite[Theorem~5]{Bou2017} we have, for every $\varepsilon>0$, 
\begin{equation}
\sum_{y=0}^{Y-1}
e\left(\frac{a(n_0+Ay)^k}{q}\right)
\ll_{D,\varepsilon}
Y^{1+\varepsilon}
\left(
q_D^{-1}+Y^{-1}+q_DY^{-D}
\right)^{1/[D(D-1)]}.
\label{Bourgain}
\end{equation}
By \eqref{qDbound} and \eqref{Ylower},
\[
q_D^{-1}\leq X^{-1/12},
\qquad
Y^{-1}\ll X^{-2/3},
\qquad
q_DY^{-D}
\ll_D X^{1/7-2D/3}.
\]
As $D \ge 5$, it follows that
\[
q_D^{-1}+Y^{-1}+q_DY^{-D}
\ll X^{-1/12}.
\]
Substituting this into \eqref{Bourgain}, we conclude that
\[
\sum_{y=0}^{Y-1}
e\left(\frac{a(n_0+Ay)^k}{q}\right)
\ll_{k,\varepsilon}
Y^{1+\varepsilon}
X^{-1/[12D(D-1)]}.
\]
Recalling \eqref{classes} and summing over $n_0$ yields
\[
\sum_{\substack{n\le X \\ (n,m)=1}} e(an^k/q) \ll X^{\eps - 1/[12D(D-1)]} \sum_{\substack{n_0 \le A \\ (n_0,m)=1}} Y_{n_0} \ll X^{1+\eps - 1/[12D(D-1)]}. 
\]

Since $U\leq M$ and $M\geq1$, we have
$
D=\lceil4U\rceil
\leq4M+1\leq5M
$
and in particular
$12D(D-1) \le 300M^2$.
Thus,
\[
\sum_{\substack{n\le X \\ (n,m)=1}} e(an^k/q)  \ll_{k,\varepsilon}
X^{1+\varepsilon-1/(300M^2)}.
\]
Choosing
\[
\varepsilon=\frac{1}{2000L^2}
\leq\frac{1}{2000M^2}
\]
furnishes \eqref{weyl}.
\end{proof}

Our main proposition is as follows.

\begin{prop}
\label{MainProp}
Let $P^\gamma\le X\le P$ and $1\le M\le L$. Let
$m$ be a square-free product of at most $R$ many
primes, all exceeding $k$, such that
\[
m\le X^{M/k}.
\]
Let $c>0$ be as in Lemma \ref{RationalFrequency}. Suppose that
$W\ge 1$ and
\[
W^2\le X^{c/(2M^2)}.
\]
Let $(c_n)_{n\le X}$ be complex coefficients
satisfying
\[
\sum_{n\le X}|c_n|^2\le X.
\]
Then, for every fixed $\beta\in\mathbb T$,
there are at most $O_k(W^2)$ many
$\alpha\in\mathbb T$ such that
$m^k \alp = \bet$ and
\[
\left|
\sum_{\substack{n\le X\\(n,m)=1}}
c_n e(\alpha n^k)
\right|
\ge\frac{X}{2W}.
\]
%The implied constant is independent of $m$, $\beta$, $W$, and the coefficients $c_n$.
\end{prop}

\begin{proof}
Fix $\beta\in\mathbb T$. Let
\[
\alpha_j=\frac{\beta+j}{m^k}
\qquad (0\le j < m^k)
\]
be the solutions $\alp \in \bT$ to
\[
m^k \alpha = \beta.
\]
Let $\cJ$ be the set of indices $j$ for which
\[
\left|
\sum_{\substack{n\le X\\(n,m)=1}}
c_n e(\alpha_j n^k)
\right|
\ge \frac{X}{2W}.
\]
We wish to show that $J := \# \cJ \ll_k W^2$.

Put
\[
\mathcal U=\{n\in\mathbb N:n\le X,\ (n,m)=1\}.
\]
Consider the matrix
\[
V=(e(\alpha_j n^k))_{\substack{j\in \cJ\\n\in\mathcal U}}.
\]
Let $K=VV^*$ be the Gram matrix of its rows, where
$V^*$ denotes the conjugate transpose of $V$. This is a Hermitian matrix of dimension $J$, with entries
\begin{equation*}
\label{eq:gram-entries}
K_{ij}
=\sum_{n\in\mathcal U}
  e((\alpha_i-\alpha_j)n^k)
= \sum_{\substack{n\le X\\(n,m)=1}}
e\left(\frac{(i-j)n^k}{m^k}\right).
\end{equation*}
Write
\[
\frac{i-j}{m^k}=\frac{a}{q},
\qquad (a,q)=1,
\]
where
\begin{equation}
\label{qdef}
q=\frac{m^k}{\gcd(i-j,m^k)}
\end{equation}
is the denominator corresponding to the entry $K_{ij}$.
Thus, each $K_{ij}$ is an exponential
sum of the type considered in Lemma \ref{RationalFrequency}.

Let $d \mid m^k$. Then $q \mid d$ if and only if $i \equiv j \mmod m^k/d$. Let $K^{(d)}$ be the matrix with entries given by
\[
K^{(d)}_{ij} = \begin{cases}
K_{ij}, &\text{if } i \equiv j \mmod m^k/d \\
0, &\text{otherwise}.
\end{cases}
\]
This is a matrix with $m^k/d$ many square blocks, namely one for each residue class. The dimension of a block is the number of $j \in \cJ$ lying in its residue class; this is at most $d$. Each block is a principal submatrix of $\tilde K^{(d)}$, which has entries
\[
K_{ij}=\sum_{\substack{n\le X\\(n,m)=1}}
e((i-j)n^k/d) \qquad 
(0\le i,j<d).
\]

The matrix $\tilde K^{(d)}$ is circulant and hence diagonalised by the discrete Fourier transform \cite[Theorem 3.2.2]{Dav1979}. Indeed, for \(b\bmod d\), the vector 
\( \bv^{(b)}=(e(jb/d))_{0\le j<d}\) is an eigenvector, since
\[
\begin{aligned}
(\tilde K^{(d)} \bv^{(b)})_i
&=\sum_{\substack{n\le X\\(n,m)=1}}e(in^k/d)
\sum_{j\bmod d}e(j(b-n^k)/d)\\
&=d\,e(ib/d)
\#\{n\le X:(n,m)=1,\ n^k\equiv b \mmod d\}.
\end{aligned}
\]
These form an orthogonal basis, and the eigenvalues of $\tilde K^{(d)}$ are precisely
\[
\lam_b := d\,\#\{n\le X:(n,m)=1,\ n^k\equiv b\mmod d\},
\]
where $b$ runs over residue classes modulo $d$.

As $d\mid m^k$, every prime divisor of $d$
exceeds $k$. For each invertible $b\bmod d$,
Hensel's lemma and the Chinese remainder
theorem tell us that there are at most
\[
k^{\omega(d)} \le k^R
\]
many solutions to $n^k\equiv b\mmod d$.
For non-invertible $b$, there are no solutions among
integers coprime to $m$. Consequently,
\[
\lambda_b\le
d k^R\left(\frac Xd+1\right)
=k^R(X+d).
\]
Any
operator-norm bound for $\tilde K^{(d)}$ also holds for its principal
sub-matrices, so
\[
\|K^{(d)} \|_{\mathrm{op}}
\le k^R(X+d).
\]
In particular,
\[
\| K^{(d)} \|_{\mathrm{op}} \ll_k X
\qquad (d\le X).
\]

Note that $K = \sum_{q \mid m^k} K^{(=q)}$, where $K^{(=q)}_{ij}$ is $K_{ij}$ if \eqref{qdef} holds and $0$ otherwise. By M\"obius inversion,
\[
K^{(=q)} = \sum_{d \mid q} \mu(q/d) K^{(d)}.
\]
As $m^k$ has at most $(k+1)^R$ many divisors,
it follows that
\[
\left \| \sum_{\substack{q \mid m^k \\ q \le X}} K^{(=q)} \right \|_{\mathrm{op}} \ll X.
\]

By Lemma \ref{RationalFrequency},
the $L^1$ norm of any row of
\[
\sum_{\substack{ q \mid m^k \\ q > X}} K^{(=q)}
\]
is $O_k(J X^{1-c/M^2})$. The matrix is Hermitian, so this is also a bound for its operator norm. With $\bc = (c_n)$,
\[
\frac{JX^2}{4W^2} \le \| V \bc \|_2^2
\le X \| V \|_2^2
= X \| K \|_{\mathrm{op}} \ll X^2 ( 1 + J X^{-c/M^2}).
\]
For $P$ large, our bound on $W^2$ now precludes the possibility that $1 < JX^{-c/M^2}$. We must instead have $1 \ge JX^{-c/M^2}$, and the desired conclusion $J \ll W^2$ emerges.
\end{proof}

\section{Main iterator}

In the context of \eqref{state} define, for $W > 0$,
\[
D_S(X,W) = X^k
\meas \{ \alp \in \bT: 
|F_S(\alp; X)| \ge X/W \}.
\]
Note that if $W < 1$ then $D_S(X, W) = 0$. Note also that $D_S(X,W)$ depends on the choice of primes indexed by $I \setminus S$. Put $K_0 = 2^{R+1}$. The following inequality is important, and is conceptually similar to \cite[Lemma 2.1]{Vau1989}. 

\begin{lem}
\label{recursion}
Let $W \ge 1$ and $\emptyset \ne T \subseteq S$. Assume that Proposition \ref{MainProp} applies to $X, W$ for some $M \in [1,L]$ and all $m \in \cD_T$. Then
\begin{align}
\notag
D_S(X,W) &\ll_k W^3 \sum_{m \in \cD_T} m^{-1} D_{S/T}(X/m, 2W) 
\\&\qquad
+
\sum_{\emptyset \ne T_1 \subseteq T} \ \sum_{m_1 \in \cD_{T_1}} m_1^k D_{S \setminus T_1} (X/m_1, K_0 W / m_1).
\label{RecursiveIneq}
\end{align}
\end{lem}

\begin{proof} Let
\[
E_S = \{ \alp \in \bT: |F_S(\alp; X)| \ge X/W \},
\]
and for $m \in \cD_T$ let
\[
H_m(\alp) = \sum_{\substack{r \in \cD_S \\ (r,m) = 1}} 
f_{X/r}(r^k \alp).
\]
We start by treating $\alp \in E_S$ such that $|H_m(\alp)| < X/(2W)$ for some $m \in \cD_T$. In this case,
\[
|F_S(\alp; X) - H_m(\alp)| \ge X/(2W).
\]
By the inclusion--exclusion principle,
\[
F_S(\alp; X) - H_m(\alp) = \sum_{\substack{r \in \cD_S \\ (r,m) > 1}} f_{X/r}(r^k \alp)
\]
is a $\pm 1$ linear combination of sub-sums over $1 < m_1 \mid m$. It follows that, for some $m_1 \in \cD_{T_1}$ with $T_1 \subseteq T$,
\[
|F_{S \setminus T_1}(m_1^k\alp; X/m_1)| \ge X / (K_0W),
\]
wherein we have used Lemma \ref{combinatorics}. Recalling that multiplication by a positive integer is measure-preserving on $\bT$, and accounting for the normalising factor in the definition of $D_S(X, W)$ we see that this case contributes the second term on the right side of \eqref{RecursiveIneq}.

Let us turn our attention to $\alp \in E_S$ such that $|H_m(\alp)| \ge X/(2W)$ for all $m \in \cD_T$, denoting the set of such $\alp$ by $E$. Write
\[
u_m(\alp) = 
\frac{m}{X} |F_{S \setminus T} (m^k \alp; X/m)| \in [0,1].
\]
By \eqref{subset}, if $\alp \in E$ then
\[
\sum_m \frac{u_m(\alp)}{m} \ge \frac 1 W,
\qquad
\sum_{\substack{m \in \cD_T \\
u_m(\alp) \ge 1/(2W)}} \frac 1m \ge \frac1{2W},
\]
whence,
\[
\sum_{m \in \cD_T}
\frac{\meas \{ \alp \in E: u_m(\alp) \ge 1/(2W) \}}{m} \ge \frac{\meas(E)}{2W}.
\]

Given $m \in \cD_T$, the quantity $u_m(\alp)$ depends on the image $\bet = m^k \alp$ rather than specifically on $\alp$. Thus, by Proposition \ref{MainProp},
\[
\meas \{ \alp \in E: u_m(\alp) \ge 1/(2W) \} \ll_k W^2 m^{-k} (X/m)^{-k} D_{S/T}(X/m, 2W).
\]
The upshot is that
\[
X^k \meas(E) \ll W^3 \sum_{m \in \cD_T} m^{-1} D_{S/T}(X/m, 2W),
\]
and so this case contributes the first term to the right hand side of \eqref{RecursiveIneq}.
\end{proof}

\section{Restricted mean value estimate}

In this section, we use Lemma \ref{recursion} to establish Theorem \ref{MeanValue} in a more explicit form. Let $c$ be as in Lemma \ref{RationalFrequency}. Fix $0< b < \min \{ 1, c/2 \}$, as well as an integer 
\[
A_0 \ge \max \{ 100/b, b^{-2} \}.
\]
For $x \in [\gam, 1]$ and $\lam \ge 0$, define
\[
B(x,\lam) = A_0 k \lam \left( 1 + \frac{\log(x/\gam)}{L} + \frac1{1 + \sqrt{\lam/x}}
\right).
\]

\begin{prop} 
\label{LargeSpectrumBound}
In the context of \eqref{state}, uniformly for $X = P^x$ and $W = P^\lam \ge 1$,
\[
D_S(X,W) \ll_k P^{4k\gam + B(x, \lam)}.
\]
\end{prop}

\begin{proof} We induct on $\# S$. The trivial estimate $D_S(X,W) \le X^k$ suffices if $x \le 4 \gam$ or $B(x,\lam) \ge kx$. By \eqref{xineqs}, this solves the base case $S = \emptyset$. We may now assume that 
\[
\# S \ge 1, \qquad
x > 4 \gam, \qquad
z := \lam/x < 1/A_0.
\]
Choose
\[
M = \begin{cases}
\min \{ L, b/ \sqrt z \}, & \text{if } z \ne 0 \\
L, &\text{if } z = 0,
\end{cases}
\]
so that $M \ge 1$ and
$2z \le 2b^2/M^2 < c/(2M^2)$. By \eqref{xineqs},
\[
\sum_{i \in S} a_i \ge x - 2 \gam > x/2.
\]
As $k$ is large and $a_i \asymp R^{-1}$ for all $i$, we can choose $T \subseteq S$ such that
\[
\frac{xM}{4k} \le \sum_{i \in T} a_i \le \frac{xM}{2k}.
\]
For $m \in \cD_T$,
\begin{equation}
\label{dbounds}
xM/(4k) \le d := \log_P m \le xM/k \le x/2,
\end{equation}
confirming the hypotheses
\[
m \le X^{M/k} \quad (m \in \cD_T), \qquad
W^2 \le X^{c/(2M^2)}
\]
of Lemma \ref{recursion}.

We claim that
\begin{equation}
\label{claim}
B(x, \lam) - B(x-d, \lam)
\ge 3 \lam.
\end{equation}
This is obvious if $\lam = 0$. To confirm it for $\lam \ne 0$, observe that for $x-d \le u \le x$ we have $\lam/u \le 2/A_0 \le 1$, whence,
\[
u \partial_u \frac1{1 + \sqrt{\lam/u}} = \frac{\sqrt{\lam/u}}{2(1+\sqrt{\lam/u})^2} \ge \frac{\sqrt z}8.
\]
Now
\begin{align*}
B(x, \lam) - B(x-d, \lam) &= \int_{x-d}^x u \partial_u B(u, \lam) \frac{\d u}{u} \ge A_0 k \lam \left( \frac 1 L + \frac{\sqrt z}{8}
\right) \int_{x-d}^x 
\frac{\d u}{u} \\
&\ge A_0 k \lam \left( \frac 1 L + \frac{\sqrt z}{8}
\right) \frac d x.
\end{align*}
Hence by \eqref{dbounds},
\[
B(x, \lam) - B(x-d, \lam)
\ge \frac{A_0 \lam M}{4}
\left( \frac 1 L + \frac{\sqrt z}{8}
\right).
\]
If $M=L$ then $B(x, \lam) - B(x-d, \lam)
\ge A_0 \lam / 4 \ge 3 \lam$. Otherwise $M \sqrt z \ge b$ and
\[
B(x, \lam) - B(x-d, \lam)
\ge A_0 \lam b / 32 \ge 3 \lam.
\]
We have confirmed the claim \eqref{claim} in all cases.

The analysis above reveals that
$
B(x,\lam) - A_0 k \lam
$
is non-decreasing in $x$, and it is readily confirmed that it is also non-decreasing in $\lam$.
Now observe that
\begin{equation}
\label{partial}
0 \le \partial_\lam B(x, \lam) \le 3 A_0 k,
\end{equation}
where the partial derivative is taken from the right at 0.

Lemma \ref{recursion} bounds $D_S(X,W) \ll_k S_1 + S_2$, where
\begin{align*}
S_1 &= W^3
\sum_{m \in \cD_T} m^{-1}
D_{S/T}(X/m, 2W), \\
S_2 &= \sum_{\emptyset \ne T_1 \subseteq T} \ \sum_{m_1 \in \cD_{T_1} } m_1^k D_{S\setminus T_1} (X/m_1, K_0 W / m_1)
\end{align*}
and $K_0 = 2^{R+1}$. By our inductive hypothesis,
\[
S_1 \ll W^3 \sum_{m \in \cD_T} m^{-1} P^{4k\gam + B(x-d, \lam + \log_P 2)}.
\]
Now by \eqref{claim}, \eqref{partial} and the mean value theorem,
\[
S_1 \ll W^3 P^{4k\gam + B(x,\lam) - 3 \lam + 3A_0 k \log_P 2} \sum_{m \in \cD_T} m^{-1} \le P^{4k\gam + B(x,\lam)} 2^{3A_0k} \ll_k P^{4k\gam + B(x,\lam)}.
\]

To bound $S_2$, write $d_1 = \log_P m_1$. For $d_1 > \lam + \log_P K_0 =: \nu$, we have 
\[
D_{S \setminus T_1}(X/m_1, K_0W/m_1) = 0.
\]
Now suppose instead that $\nu - d_1 \ge 0$. We also have $d_1 \le d \le x/2$ and hence 
\[
x-d_1 \ge x/2 > \lam.
\]
As $B(x,\lam) - A_0 k \lam$ is non-decreasing in each variable,
\[
B(x-d_1, \nu-d_1) \le B(x,\nu) - A_0 k d_1.
\]
Finally, by our inductive hypothesis,
\begin{align*}
S_2 &\ll \sum_{\emptyset \ne T_1 \subseteq T} \ \sum_{\substack{m_1 \in \cD_{T_1} \\ d_1 \le \nu}} m_1^k
P^{4k\gam + B(x-d_1, \lam + \log_P K_0 - d_1)}
\\ &\ll
P^{4k\gam + B(x,\lam)} \sum_{m_1=1}^\infty m_1^{(1-A_0)k} \ll P^{4k\gam + B(x,\lam)}.
\end{align*}
\end{proof}

\begin{cor} 
\label{MVE}
For $r= 3 A_0 k+ 1$ and $F$ as in \eqref{state},
\[
\int_\bT |F(\alp)|^r \d \alp \ll_k P^{r-k + 4/k^2}.
\]
\end{cor}

\begin{proof} For $t \in (0,1]$, an application of Proposition \ref{LargeSpectrumBound} with $W = 1/t$ delivers
\[
P^k \meas \{ \alp \in \bT: |F(\alp)| \ge tP \} \ll P^{4k\gam + B(1,\lam)},
\]
where $W = P^\lam$. As $B(1, \lam) \le 3A_0 k \lam$, we obtain
\[
\meas \{ \alp \in \bT: |F(\alp)| \ge tP \} \ll P^{(4\gam - 1)k} t^{-3A_0k}.
\]
Since
\[
(|F(\alp)|/P)^r = \int_0^{|F(\alp)|/P} rt^{r-1} \d t = \int_0^1 rt^{r-1} 1_{|F(\alp)| \ge tP} \d t,
\]
we have
\begin{align*}
\int_\bT |F(\alp)|^r \d \alp &\ll P^r \int_0^1 r t^{r-1} \meas \{ \alp \in \bT: |F(\alp)| \ge tP \}  \d t \\
&\ll P^{r+(4\gam - 1)k} \int_0^1 r t^{r-1-3A_0k} \d t.
\end{align*}
As $k\gam = k^{-2}$ and $r= 3A_0k+1$, this is $O(P^{r-k+4/k^2})$.
\end{proof}

\section{The circle method}

In this section, we complete the proof of Theorem \ref{MainThm}. For $1 \le Q \le P$, define major arcs
\[
\fM(Q) = \bigcup_{\substack{1 \le a \le q \le Q \\ (a,q)=1}} \{ \alp \in \bT: \| \alp - a/q \| \le Q/(qP^k) \},
\]
and write $\fm = \bT \setminus \fM(P)$ for the minor arcs. Let
\[
S(q,a) = \sum_{x \le q} e(ax^k / q), \qquad
V(\xi) = \int_0^1 e(\xi t^k) \d t.
\]

\begin{lem} [Smooth Weyl sum estimates]
\label{SmoothWeyl}
Fix $k \in \bN$ sufficiently large and $\eta > 0$ sufficiently small. Then:
\begin{enumerate}[(a)]
\item With $\rho_k^{-1} = k(\log k + 8.02113)$,
\[
\sup_{\alp \in \fm} |f_P(\alp)| \ll_k P^{1-\rho_k}.
\]
\item For some $\kap_k > 0$ we have, uniformly for $1 \le Q \le P$,
\[
\int_{\fM(P) \setminus \fM(Q)} |f_P(\alp)|^{2k} \d \alp \ll_k P^k Q^{-\kap_k}.
\]
\item There exists $c_\eta > 0$ such that for any fixed $Q$ we have, uniformly for
\begin{equation}
\label{HeightQ}
q \le Q, \qquad (a,q) = 1, \qquad |\xi| \le Q,
\end{equation}
the approximation
\[
f_Y(a/q + \xi/Y^k) = c_\eta Y q^{-1} S(q,a) V(\xi) + o(Y) \qquad (Y \to \infty)
\]
\end{enumerate}
\end{lem}

\begin{proof} 
\begin{enumerate}[(a)]
\item Apply \cite[Theorem 1.1]{BW2025}.
\item Apply \cite[Theorem 1.3]{BW2024} with 
\[
u = 2k, \qquad t = \lfloor k/2 \rfloor, \qquad \ome = \ome' = 1.
\]
\item By the estimate used immediately before \cite[Equation (6.15)]{BW2023}, we have
\[
f_Y(a/q + \xi/Y^k) = c_\eta q^{-1} S(q,a) v(\xi/Y^k) + o(Y), 
\]
where $v(\bet) = k^{-1} \sum_{m \le Y^k} m^{-1 + 1/k} e(\bet m)$. As explained on Page 16 of \cite{Vaughan97},
\[
v(\bet) + O(1 + Y^k|\bet|) = \int_0^Y e(\bet \gam^k) \d \gam = Y V(\bet Y^k).
\]
\end{enumerate}
\end{proof}

Recall the small, positive constant $\ome$ from \eqref{omega}, and let $\tau = (\log(1 + \ome))^R$. Also recall that
$
\cA_k(P) = \{ mu \le P: m \in \cD_I, \ u \text{ is } (P/m)^\eta\text{-smooth} \}.
$

\begin{lem}
[Major arc asymptotic]
\label{major}
For any fixed $Q$ we have, uniformly for \eqref{HeightQ},
\[
F(a/q + \xi / P^k) =
c_\eta \tau P q^{-1} S(q,a) V( \xi) + o(P) \qquad (P \to \infty).
\]
In particular, 
\[
\# \cA_k(P) = (c_\eta \tau + o(1)) P.
\]
\end{lem}

\begin{proof}
Recalling the setup from \S \ref{setup}, we have $\sum a_i \le 1-2\gam$.
For $m \in \cD_I$, we apply Part (c) of Lemma \ref{SmoothWeyl} with $Y = P/m \ge P^\gam$. Observe that $(m,q) = 1$, since $P$ is large, so $S(q, m^k a) = S(q,a)$ by permuting the residue classes. Summing over $m$ yields
\begin{align*}
F(a/q + \xi/P^k) &= \sum_{m \in \cD_I} f_Y(m^k a/q + \xi/ Y^k)
\\ &= \sum_{m \in \cD_I} (c_\eta P m^{-1} q^{-1} S(q,a) V(\xi) + o(P/m)).
\end{align*}
It follows from Mertens's theorem that 
\[
\sum_{m \in \cD_I} m^{-1} \sim \tau \qquad (P \to \infty),
\]
and substituting this completes the proof.
\end{proof}

Note that Corollary \ref{MVE} and Lemma \ref{major} imply Theorem \ref{MeanValue} with $d_k = c_\eta \tau$.

\begin{proof}
[Proof of Theorem \ref{MainThm}]
Let $N$ be a large, positive integer, and put
\[
P = N^{1/k}, \qquad
r = 3A_0k + 1, \qquad
s = 2k + r.
\]
By orthogonality,
\[
R_{s,k}(N) \ge \int_\bT f_P(\alp)^{2k} F(\alp)^r e(-N\alp) \d \alp =: \tilde R(N).
\]
By Part (a) of Lemma \ref{SmoothWeyl}, and Corollary \ref{MVE},
\begin{align}
\notag
\int_\fm |f_P(\alp)^{2k} F(\alp)^r| \d \alp &\le (\sup_{\alp \in \fm} |f_P(\alp)|)^{2k} \int_\bT |F(\alp)|^r \d \alp \\
&\ll_k P^{s-k+4/k^2 - 2/ (\log k + 8.02113)} = o(P^{s-k}).
\label{minor}
\end{align}
The trivial bound $|F(\alp)| \le P$, together with Part (b) of Lemma \ref{SmoothWeyl}, yields
\begin{equation}
\label{prune}
\int_{\fM(P) \setminus \fM(Q)} |f_P(\alp)^{2k} F(\alp)^r| \ll P^{s-k} Q^{-\kap_k},
\end{equation}
uniformly for $1 \le Q \le P$.

We take $Q$ to be large but fixed, and define an enlarged set of major arcs
\[
\fK(Q) = \bigcup_{\substack{1 \le a \le q \le Q \\ (a,q)=1}} \left \{
\alp \in \bT: \left \|
\alp - \frac a q
\right \| \le Q/P^k
\right \}.
\]
Define
\[
\fS_{s,Q}(N) = \sum_{q \le Q} \ \sum_{\substack{a \mmod q \\ (a,q) = 1}} (q^{-1} S(q,a))^s e(-Na/q), \qquad
\fJ_{s,Q} = \int_{-Q}^Q V(\xi)^s e(-\xi) \d \xi.
\]
Using Part (c) of Lemma \ref{SmoothWeyl} for $f_P$, and using Lemma \ref{major} for $F$, we obtain
\[
\int_{\fK(Q)} f_P(\alp)^{2k} F(\alp)^r e(-N\alp) \d \alp
= c_\eta^s \tau^r P^s \fS_{s,Q}(N) \fJ'_{s,Q} + o(P^{s-k}),
\]
where
\[
\fJ'_{s,Q} = \int_{|\bet| \le QP^{-k}} V(P^k \bet)^s e(-N\bet)
\d \bet = P^{-k} \fJ_{s,Q}.
\]
Thus,
\begin{equation}
\label{MajorContribution}
\int_{\fK(Q)} f_P(\alp)^{2k} F(\alp)^r
e(-N\alp) \d \alp
= c_\eta^s \tau^r \fS_{s,Q}(N) \fJ_{s,Q}
P^{s-k} + o(P^{s-k}).
\end{equation}

Combining \eqref{minor}, \eqref{prune} and \eqref{MajorContribution} furnishes
\begin{equation*}
\label{conc}
R_{s,k}(N) \ge \tilde R(N) = c_\eta^s \tau^r \fS_{s,Q}(N) \fJ_{s,Q}
P^{s-k} + o(P^{s-k}) + O(P^{s-k} Q^{-\kap_k}).
\end{equation*}
The classical theory \cite[Chapter 4]{Vaughan97} tells us that $\fS_{s,Q}(N), \fJ_{s,Q}$ converge uniformly in $N$ to $\fS_s(N), \fJ_s > 0$ respectively, and that $\fS_s(N) \gg_s 1$.
The implied constant in $O(P^{s-k} Q^{-\kap_k})$ does not depend on $Q$. Taking $Q$ to be a sufficiently large constant, and then taking $P$ large, we conclude that 
$
R_{s,k}(N) \gg P^{s-k}.
$
\end{proof}

\providecommand{\bysame}{\leavevmode\hbox to3em{\hrulefill}\thinspace}

\end{document}